\documentclass[reqno, 12pt]{amsart}
\usepackage{a4wide}
\usepackage{amscd}
\input xy
\xyoption{all}

\newtheorem{thm}{Theorem}[section]

\newtheorem{prop}[thm]{Proposition}
\newtheorem{cor}[thm]{Corollary}
\newtheorem{lm}[thm]{Lemma}

\numberwithin{equation}{section}

\begin{document}

\title{Real hypersurfaces with isometric Reeb flow in complex quadrics }
\author{\textsc{J\"{u}rgen Berndt and Young Jin Suh}}

\address{King's College London\\ Department of Mathematics \\  London WC2R 2LS \\ United Kingdom}
\address{Kyungpook National University \\ College of Natural Sciences\\  Department of Mathematics \\
Daegu 702-701\\ Republic of Korea}
\date{}

\begin{abstract}
We classify real hypersurfaces with isometric Reeb flow in the complex quadrics $Q^m = SO_{m+2}/SO_mSO_2$, $m \geq 3$. We show that $m$ is even, say $m = 2k$, and any such hypersurface is an open part of a tube around a $k$-dimensional complex projective space ${\mathbb C}P^k$ which is embedded canonically in $Q^{2k}$ as a totally geodesic complex submanifold. As a consequence we get the non-existence of real hypersurfaces with isometric Reeb flow in odd-dimensional complex quadrics $Q^{2k+1}$, $k \geq 1$. To our knowledge the odd-dimensional complex quadrics are the first examples of homogeneous K\"{a}hler manifolds which do not admit a real hypersurface with isometric Reeb flow.
\end{abstract}

\maketitle
\thispagestyle{empty}

\footnote[0]{This work was supported by grant Proj. No. NRF-2011-220-C00002 from National Research Foundation of Korea and the second author by
grant Proj. No. BSRP-2012-R1A2A2A-01043023.
}

\section{Introduction}
Let $M$ be a real hypersurface in a K\"{a}hler manifold $\bar{M}$. The complex structure $J$ on $\bar{M}$ induces locally an almost contact metric structure $(\phi,\xi,\eta,g)$ on $M$. In the context of contact geometry, the unit vector field $\xi$ is often referred to as the Reeb vector field on $M$ and its flow is known as the Reeb flow. The Reeb flow has been of significant interest in recent years, for example in relation to the Weinstein Conjecture. We are interested in the Reeb flow in the context of Riemannian geometry, namely in the classification of real hypersurfaces with isometric Reeb flow in homogeneous K\"{a}hler manifolds. 

For the complex projective space ${\mathbb C}P^m$ a full classification was obtained by Okumura in \cite{O}. He proved that the Reeb flow on a real hypersurface in ${\mathbb C}P^m = SU_{m+1}/S(U_mU_1)$ is isometric if and only if $M$ is an open part of a tube around a totally geodesic ${\mathbb C}P^k \subset {\mathbb C}P^m$ for some $k \in \{0,\ldots,m-1\}$. For the complex $2$-plane Grassmannian $G_2({\mathbb C}^{m+2})= SU_{m+2}/S(U_mU_2)$ the classification was obtained by the authors in \cite{BS}. The Reeb flow on a real hypersurface in $G_2({\mathbb C}^{m+2})$ is isometric if and only if $M$ is an open part of a tube around a totally geodesic $G_2({\mathbb C}^{m+1}) \subset G_2({\mathbb C}^{m+2})$. In this paper we investigate this problem for the complex quadric $Q^m = SO_{m+2}/SO_mSO_2$. In view of the previous two results a natural expectation is that the classification involves at least the totally geodesic $Q^{m-1} \subset Q^m$. Surprisingly, this is not the case. Our main result states:

\begin{thm} \label{mainresult}
Let $M$ be a real hypersurface of the complex quadric $Q^m$, $m\geq 3$. The Reeb flow on $M$ is isometric if and only if $m$ is even, say $m = 2k$, and $M$ is an open part of a tube around a totally geodesic ${\mathbb C}P^k \subset Q^{2k}$. 
\end{thm}

Every tube around a totally geodesic ${\mathbb C}P^k \subset Q^{2k}$ is a homogeneous hypersurface. In fact, the closed subgroup $U_{k+1}$ of $SO_{2k+2}$ acts on $Q^{2k}$ with cohomogeneity one. The two singular orbits are totally geodesic ${\mathbb C}P^k \subset Q^{2k}$ and the principal orbits are the tubes around any of these two singular orbits. So as a corollary we get:

\begin{cor}
Let $M$ be a connected complete real hypersurface in the complex quadric $Q^{2k}$, $k \geq 2$. If the Reeb flow on $M$ is isometric, then $M$ is a homogeneous hypersurface of $Q^{2k}$.
\end{cor}

It is remarkable that in this situation the existence of a particular one-parameter group of isometries implies transitivity of the isometry group. As another interesting consequence we get:

\begin{cor} 
There are no real hypersurfaces with isometric Reeb flow in the odd-dimensional complex quadric $Q^{2k+1}$, $k \geq 1$.
\end{cor}

To our knowledge the odd-dimensional complex quadrics are the first examples of homogeneous K\"{a}hler manifolds which do not admit a real hypersurface with isometric Reeb flow.

Our paper is organized as follows. In Section 2 we present basic material about the complex quadric $Q^m$, including its Riemannian curvature tensor and a description of its singular tangent vectors. Apart from the complex structure $J$
there is another distinguished geometric structure on $Q^m$, namely a parallel rank two vector bundle ${\mathfrak A}$ which contains an $S^1$-bundle of real structures on the tangent spaces of $Q^m$. This geometric structure determines a maximal ${\mathfrak A}$-invariant subbundle ${\mathcal Q}$ of the tangent bundle $TM$ of a real hypersurface $M$ in $Q^m$. In Section 3 we investigate the geometry of this subbundle ${\mathcal Q}$. In Section 4 we describe the canonical embedding of ${\mathbb C}P^k$ into $Q^{2k}$ as a totally geodesic complex submanifold and investigate the geometry of the tubes around ${\mathbb C}P^k$. We will show that the Reeb flow on these tubes is isometric. In Section 5 we determine some geometric consequences of the Codazzi equation for real hypersurfaces in $Q^m$. 

Finally, in Section 6, we present the proof of Theorem \ref{mainresult}. The first step is to prove that the normal bundle of a real hypersurface $M$ with isometric Reeb flow consists of a particular type of singular tangent vectors of $Q^m$, the so-called ${\mathfrak A}$-isotropic vectors. The next step is to show that ${\mathcal Q}$ is invariant under the shape operator of $M$. Putting all this information into the Codazzi equation then allows us to compute explicitly the principal curvatures and principal curvatures spaces of $M$. A particular consequence of this is that $m$ is even, say $m = 2k$.
Using Jacobi field theory we then show that $M$ has a smooth focal manifold at a constant distance which is embedded in $Q^{2k}$ as a totally geodesic complex submanifold of complex dimension $k$. Using Klein's \cite{K} classification of totally geodesic submanifolds in complex quadrics we will then show that this focal manifold is a totally geodesic ${\mathbb C}P^k \subset Q^{2k}$.

\section{The complex quadric} \label{quadric}

The homogeneous quadratic equation $z_1^2 + \ldots + z_{m+2}^2 = 0$
on ${\mathbb C}^{m+2}$ defines a complex hypersurface $Q^m$ in the $(m+1)$-dimensional complex projective space ${\mathbb C}P^{m+1} = SU_{m+2}/S(U_{m+1}U_1)$. The hypersurface $Q^m$ is known as the $m$-dimensional complex quadric. The complex structure $J$ on ${\mathbb C}P^{m+1}$ naturally induces a complex structure on  $Q^m$ which we will denote by $J$ as well. We equip $Q^m$ with the Riemannian metric $g$ which is induced from the Fubini Study metric on ${\mathbb C}P^{m+1}$ with constant holomorphic sectional curvature $4$. The $1$-dimensional quadric $Q^1$ is isometric to the round $2$-sphere $S^2$. For $m \geq 2$ the triple $(Q^m,J,g)$ is a Hermitian symmetric space of rank two and its maximal sectional curvature is equal to $4$. The $2$-dimensional quadric $Q^2$ is isometric to the Riemannian product $S^2 \times S^2$. We will assume $m \geq 3$ for the main part of this paper.

For a nonzero vector $z \in {\mathbb C}^{m+1}$ we denote by $[z]$ the complex span of $z$, that is, $[z] = \{\lambda z \mid \lambda \in {\mathbb C}\} $.
Note that by definition $[z]$ is a point in ${\mathbb C}P^{m+1}$.
As usual, for each $[z] \in {\mathbb C}P^{m+1}$ we identify $T_{[z]}{\mathbb C}P^{m+1}$ with the orthogonal complement ${\mathbb C}^{m+2} \ominus [z]$ of $[z]$ in ${\mathbb C}^{m+2}$. For $[z] \in Q^m$ the tangent space $T_{[z]}Q^m$ can then be identified canonically with the orthogonal complement ${\mathbb C}^{m+2} \ominus ([z] \oplus [\bar{z}])$ of $[z] \oplus [\bar{z}]$ in ${\mathbb C}^{m+2}$. Note that $\bar{z} \in \nu_{[z]}Q^m$ is a unit normal vector of $Q^m$ in ${\mathbb C}P^{m+1}$ at the point $[z]$.

We denote by $A_{\bar{z}}$ the shape operator of $Q^m$ in ${\mathbb C}P^{m+1}$ with respect to $\bar{z}$. Then we have $A_{\bar{z}}w = \overline{w}$ for all $w \in T_{[z]}Q^m$, that is,  $A_{\bar{z}}$ is just complex conjugation restricted to $T_{[z]}Q^m$.
The shape operator $A_{\bar{z}}$ is an antilinear involution on the complex vector space $T_{[z]}Q^m$ and
$$T_{[z]}Q^m = V(A_{\bar{z}}) \oplus JV(A_{\bar{z}}),$$
where $V(A_{\bar{z}}) =  {\mathbb R}^{m+2} \cap T_{[z]}Q^m$ is the $(+1)$-eigenspace and $JV(A_{\bar{z}}) = i{\mathbb R}^{m+2} \cap T_{[z]}Q^m$ is the $(-1)$-eigenspace of $A_{\bar{z}}$.  Geometrically this means that the shape operator $A_{\bar{z}}$ defines a real structure on the complex vector space $T_{[z]}Q^m$. Recall that a real structure on a complex vector space $V$ is by definition an antilinear involution $A : V \to V$. Since the normal space $\nu_{[z]} Q^m$ of $Q^m$ in ${\mathbb C}P^{m+1}$ at ${[z]}$ is a complex subspace of $T_{[z]}{\mathbb C}P^{m+1}$ of complex dimension one, every normal vector in $\nu_{[z]} Q^m$ can be written as $\lambda \bar{z}$ with some $\lambda \in  {\mathbb C}$. The shape operators $A_{\lambda\bar{z}}$ of $Q^m$ define a rank two vector subbundle ${\mathfrak A}$ of the endomorphism bundle ${\rm End}(TQ^m)$. Since the second fundamental form of the embedding $Q^m \subset {\mathbb C}P^{m+1}$ is parallel (see e.g.\ \cite{S}), ${\mathfrak A}$ is a parallel subbundle of ${\rm End}(TQ^m)$. For $\lambda \in S^1 \subset {\mathbb C}$ we again get a real structure $A_{\lambda\bar{z}}$ on $T_{[z]}Q^m$ and we have $V(A_{\lambda\bar{z}}) = \lambda V(A_{\bar{z}})$. We thus have an $S^1$-subbundle of ${\mathfrak A}$ consisting of real structures on the tangent spaces of $Q^m$.

The Gauss equation for the complex hypersurface $Q^m \subset {\mathbb C}P^{m+1}$ implies that the Riemannian curvature tensor $R$ of $Q^m$ can be expressed in terms of the Riemannian metric $g$, the complex structure $J$ and a generic real structure $A$ in ${\mathfrak A}$:
\begin{eqnarray*}
R(X,Y)Z & = & g(Y,Z)X - g(X,Z)Y + g(JY,Z)JX - g(JX,Z)JY - 2g(JX,Y)JZ \\
 & & +\, g(AY,Z)AX - g(AX,Z)AY + g(JAY,Z)JAX - g(JAX,Z)JAY.
\end{eqnarray*}
Note that the complex structure $J$ anti-commutes with each endomorphism $A \in {\mathfrak A}$, that is, $AJ = -JA$.

A nonzero tangent vector $W \in T_{[z]}Q^m$ is called singular if it is tangent to more than one maximal flat in $Q^m$. There are two types of singular tangent vectors for the complex quadric $Q^m$:
\begin{itemize}
\item[1.] If there exists a real structure $A \in {\mathfrak A}_{[z]}$ such that $W \in V(A)$, then $W$ is singular. Such a singular tangent vector is called ${\mathfrak A}$-principal.
\item[2.] If there exist a real structure $A \in {\mathfrak A}_{[z]}$ and orthonormal vectors $X,Y \in V(A)$ such that $W/||W|| = (X+JY)/\sqrt{2}$, then $W$ is singular. Such a singular tangent vector is called ${\mathfrak A}$-isotropic.
\end{itemize}
Basic complex linear algebra shows that for every unit tangent vector $W \in T_{[z]}Q^m$ there exist a real structure $A \in {\mathfrak A}_{[z]}$ and orthonormal vectors $X,Y \in V(A)$ such that
\[
W = \cos(t)X + \sin(t)JY
\]
for some $t \in [0,\pi/4]$. The singular tangent vectors correspond to the values $t = 0$ and $t = \pi/4$.

\section{The maximal ${\mathfrak A}$-invariant subbundle ${\mathcal Q}$ of $TM$}

Let $M$ be a  real hypersurface in $Q^m$ and denote by $(\phi,\xi,\eta,g)$ the induced almost contact metric structure on $M$ and by $\nabla$ the induced Riemannian connection on $M$. Note that $\xi = -JN$, where $N$ is a (local) unit normal vector field of $M$. The vector field $\xi$ is known as the Reeb vector field of $M$. If the integral curves of $\xi$ are geodesics in $M$, the hypersurface $M$ is called a Hopf hypersurface. The integral curves of $\xi$ are geodesics in $M$ if and only if $\xi$ is a principal curvature vector of $M$ everywhere.
The tangent bundle $TM$ of $M$ splits orthogonally into  $TM = {\mathcal C} \oplus {\mathcal F}$, where ${\mathcal C} = {\rm ker}(\eta)$ is the maximal complex subbundle of $TM$ and ${\mathcal F} = {\mathbb R}\xi$. The structure tensor field $\phi$ restricted to ${\mathcal C}$ coincides with the complex structure $J$ restricted to ${\mathcal C}$, and we have $\phi \xi = 0$. We denote by $\nu M$ the normal bundle of $M$.

We first introduce some notations. For a fixed real structure $A \in {\mathfrak A}_{[z]}$ and $X \in T_{[z]}M$ we decompose $AX$ into its tangential and normal component, that is,
\[
AX = BX + \rho(X)N
\]
where $BX$ is the tangential component of $AX$ and 
\[
\rho(X) = g(AX,N) = g(X,AN) = g(X,AJ\xi) = -g(X,JA\xi) = g(JX,A\xi).
\]
Since $JX = \phi X + \eta(X)N$ and $A\xi = B\xi + \rho(\xi)N$ we also have
\[
\rho(X) = g(\phi X,B\xi) + \eta(X)\rho(\xi)  = \eta(B\phi X) + \eta(X)\rho(\xi).
\]
We also define 
\[
\delta = g(N,AN) = g(JN,JAN) = -g(JN,AJN) = -g(\xi,A\xi).
\]

At each point $[z] \in M$ we define
\[
{\mathcal Q}_{[z]} = \{X \in T_{[z]}M \mid AX \in T_{[z]}M\ {\rm for\ all}\ A \in {\mathfrak A}_{[z]}\},
\]
which is the maximal ${\mathfrak A}_{[z]}$-invariant subspace of $T_{[z]}M$. 

\begin{lm} \label{principal}
The following statements are equivalent:
\begin{itemize}
\item[(i)] The normal vector $N_{[z]}$ is ${\mathfrak A}$-principal,
\item[(ii)] ${\mathcal Q}_{[z]} = {\mathcal C}_{[z]}$,
\item[(iii)] There exists a real structure $A \in {\mathfrak A}_{[z]}$ such that $AN_{[z]} \in {\mathbb C}\nu_{[z]}M$.
\end{itemize}
\end{lm}

\proof We first assume that $N_{[z]}$ is ${\mathfrak A}$-principal. By definition, this means that there exists a real structure $A \in {\mathfrak A}_{[z]}$ such that $N_{[z]} \in V(A)$, that is, $AN_{[z]} = N_{[z]}$. Then we have $A\xi_{[z]} = -AJN_{[z]} = JAN_{[z]} = JN_{[z]} = -\xi_{[z]}$. It follows that $A$ restricted to the plane ${\mathbb C}\nu_{[z]}M$ is the orthogonal reflection in the line $\nu_{[z]}M$. Since all endomorphisms in ${\mathfrak A}_{[z]}$ differ just by scalar multiplication with a complex number we see that ${\mathbb C}\nu_{[z]}M$ is invariant under ${\mathfrak A}_{[z]}$. This implies that ${\mathcal C}_{[z]} = T_{[z]}Q^m \ominus {\mathbb C}\nu_{[z]}M$ is invariant under ${\mathfrak A}_{[z]}$, and hence ${\mathcal Q}_{[z]} = {\mathcal C}_{[z]}$. 

Conversely, assume that ${\mathcal Q}_{[z]} = {\mathcal C}_{[z]}$. Then ${\mathbb C}\nu_{[z]}M = T_{[z]}Q^m \ominus {\mathcal C}_{[z]}$ is invariant under ${\mathfrak A}_{[z]}$. Since the real dimension of ${\mathbb C}\nu_{[z]}M$ is $2$ and all real structures in ${\mathfrak A}_{[z]}$ differ just by scalar multiplication with a unit complex number, there exists a real structure $A \in {\mathfrak A}_{[z]}$ which fixes $N_{[z]}$. Thus $N_{[z]}$ is ${\mathfrak A}$-principal. This proves the equivalence of (i) and (ii). 

The equivalence of (ii) and (iii) follows from the fact that all real structures in ${\mathfrak A}_{[z]}$ differ just by scalar multiplication with a unit complex number.
\hfill $\Box$

\medskip
Assume now that the normal vector $N_{[z]}$ is not ${\mathfrak A}$-principal.
Then there exists a real structure $A \in {\mathfrak A}_{[z]}$ such that
\[ N_{[z]} = \cos(t)Z_1 + \sin(t)JZ_2 \]
for some orthonormal vectors $Z_1,Z_2 \in V(A)$ and $0 < t \leq \frac{\pi}{4}$. This implies
\begin{eqnarray*}
N_{[z]} & = & \cos(t)Z_1 + \sin(t)JZ_2, \\
AN_{[z]} & = & \cos(t)Z_1 - \sin(t)JZ_2, \\
\xi_{[z]} & = & \sin(t)Z_2 - \cos(t)JZ_1, \\
A\xi_{[z]} & = & \sin(t)Z_2 + \cos(t)JZ_1,
\end{eqnarray*}
and therefore ${\mathcal Q}_{[z]} = T_{[z]}Q^m \ominus ([Z_1] \oplus [Z_2])$ is strictly contained in ${\mathcal C}_{[z]}$. Moreover, we have
\[
A\xi_{[z]} = B\xi_{[z]}\ {\rm and}\  \rho(\xi_{[z]}) = 0.
\]
We have
\[
g(B\xi_{[z]} + \delta \xi_{[z]},N_{[z]}) = 0,\ g(B\xi_{[z]} + \delta \xi_{[z]},\xi_{[z]}) = 0,\ g(B\xi_{[z]} + \delta\xi_{[z]},B\xi_{[z]} + \delta\xi_{[z]}) = \sin^2(2t),
\]
and therefore
\[ U_{[z]} = \frac{1}{\sin(2t)}(B\xi_{[z]} + \delta\xi_{[z]}) \]
is a unit vector in ${\mathcal C}_{[z]}$ and
\[
{\mathcal C}_{[z]} = {\mathcal Q}_{[z]} \oplus [U_{[z]}]\ {\rm (orthogonal\ sum)}.
\]
If $N_{[z]}$ is not ${\mathfrak A}$-principal at $[z]$, then $N$ is not ${\mathfrak A}$-principal in an open neighborhood of $[z]$, and therefore $U$ is a well-defined unit vector field on that open neighborhood. We summarize this in the following 

\begin{lm} \label{notprincipal}
Assume that $N_{[z]}$ is not ${\mathfrak A}$-principal at $[z]$. Then there exists an open neighborhood of $[z]$ in $M$ and a section $A$ in ${\mathfrak A}$ on that neighborhood consisting of real structures such that
\begin{itemize}
\item[(i)] $A\xi = B\xi$ and $\rho(\xi) = 0$,
\item[(ii)] $U = (B\xi + \delta\xi)/||B\xi + \delta\xi||$ is a unit vector field tangent to ${\mathcal C}$
\item[(iii)] ${\mathcal C} = {\mathcal Q} \oplus [U]$.
\end{itemize}
\end{lm}

In the following we will always assume that $A$ is chosen in line with Lemma \ref{notprincipal} in neighborhoods where the normal vector field $N$ is not ${\mathfrak A}$-principal.

\section {Tubes around the totally geodesic ${\mathbb C}P^k \subset Q^{2k}$}

We assume that $m$ is even, say $m = 2k$. The map
$$
{\mathbb C}P^k \to Q^{2k} \subset {\mathbb C}P^{2k+1} \ ,\  [z_1,\ldots,z_{k+1}] \mapsto [z_1,\ldots,z_{k+1},iz_1,\ldots,iz_{k+1}]
$$
provides an embedding of ${\mathbb C}P^k$ into $Q^{2k}$ as a totally geodesic complex submanifold.

We define a complex structure $j$ on ${\mathbb C}^{2k+2}$ by
$$
j(z_1,\ldots,z_{k+1},z_{k+2},\ldots,z_{2k+2}) = (-z_{k+2},\ldots,-z_{2k+2},z_1,\ldots,z_{k+1}).
$$
Note that $ij = ji$. We can then identify ${\mathbb C}^{2k+2}$ with ${\mathbb C}^{k+1} \oplus j{\mathbb C}^{k+1}$ and get
$$
T_{[z]}{\mathbb C}P^k = \{X + jiX \mid X \in {\mathbb C}^{k+1} \ominus [z]\} = \{X + ijX \mid X \in V(A_{\bar{z}})\}.
$$
Note that the complex structure $i$ on ${\mathbb C}^{2k+2}$ corresponds to the complex structure $J$ on $T_{[z]}Q^{2k}$ via the obvious identifications. For the normal space $\nu_{[z]}{\mathbb C}P^k$ of ${\mathbb C}P^k$ at $[z]$ we have
\[
\nu_{[z]}{\mathbb C}P^k = A_{\bar{z}}(T_{[z]}{\mathbb C}P^k) = \{X - ijX \mid X \in V(A_{\bar{z}})\}.
\]
It is easy to see that both the tangent bundle and the normal bundle of ${\mathbb C}P^k$ consist of ${\mathfrak A}$-isotropic singular tangent vectors of $Q^{2k}$.

We will now calculate the principal curvatures and principal curvature spaces of the tube with radius $0 < r < \pi/2$ around ${\mathbb C}P^k$ in $Q^{2k}$. Let $N$ be a unit normal vector of ${\mathbb C}P^k$ in $Q^{2k}$ at ${[z]} \in {\mathbb C}P^k$. Since $N$ is ${\mathfrak A}$-isotropic, the four vectors $N,JN,AN,JAN$ are pairwise orthonormal and the normal Jacobi operator $R_N$ is given by
\[
R_N Z = R(Z,N)N  =   Z - g(Z,N)N + 3g(Z,JN)JN   - g(Z,AN)AN - g(Z,JAN)JAN.
\]
This implies readily that $R_N$ has the three eigenvalues $0,1,4$ with corresponding eigenspaces ${\mathbb R}N\oplus [AN]$, $T_{[z]}Q^{2k} \ominus ([N] \oplus [AN])$ and ${\mathbb R}JN$. Since $[N] \subset \nu_{[z]}{\mathbb C}P^k$ and $[AN] \subset  T_{[z]}{\mathbb C}P^k$, we conclude that both $T_{[z]}{\mathbb C}P^k$ and $\nu_{[z]}{\mathbb C}P^k$ are invariant under $R_N$. 

To calculate the principal curvatures of the tube of radius $0 < r <\pi/2$ around ${\mathbb C}P^k$ we use the Jacobi field method as described in Section  8.2 of \cite{BCO}. Let $\gamma$ be the geodesic  in $Q^{2k}$ with $\gamma(0) = [z]$ and $\dot{\gamma}(0) = N$ and denote by $\gamma^\perp$ the parallel subbundle of $TQ^{2k}$ along $\gamma$ defined by $\gamma^\perp_{\gamma(t)} = T_{[\gamma(t)]}Q^{2k} \ominus {\mathbb R}\dot{\gamma}(t)$. Moreover, define the $\gamma^\perp$-valued tensor field $R^\perp_{\gamma}$ along $\gamma$ by $R^\perp_{\gamma(t)}X = R(X,\dot{\gamma}(t))\dot{\gamma}(t)$.
Now consider the ${\rm End}(\gamma^\perp)$-valued differential equation
\[
Y^{\prime\prime} + R^\perp_{\gamma} \circ Y = 0.
\]
Let $D$ be the unique solution of this differential equation with initial values
\[
D(0) = \begin{pmatrix} I & 0 \\ 0 & 0 \end{pmatrix}\ ,\ D^\prime(0) = \begin{pmatrix} 0 & 0 \\ 0 & I \end{pmatrix},
\]
where the decomposition of the matrices is with respect to
\[
\gamma^\perp_{[z]} = T_{[z]}{\mathbb C}P^k \oplus (\nu_{[z]}{\mathbb C}P^k \ominus {\mathbb R}N)
\]
and $I$ denotes the identity transformation on the corresponding space. Then the shape operator $S(r)$ of the tube of radius $0 < r <\pi/2$ around ${\mathbb C}P^k$ with respect to $-\dot{\gamma}(r)$ is given by
\[
S(r) = D^\prime(r)\circ D^{-1}(r).
\]
If we decompose $\gamma^\perp_{[z]}$ further into
\[
\gamma^\perp_{[z]} = (T_{[z]}{\mathbb C}P^k \ominus [AN]) \oplus [AN] \oplus (\nu_{[z]}{\mathbb C}P^k \ominus [N]) \oplus {\mathbb R}JN,
\]
we get by explicit computation that
\[
S(r) = \begin{pmatrix}
-\tan(r) & 0 & 0 & 0 \\
0 & 0 & 0 & 0 \\
0 & 0 & \cot(r) & 0 \\
0 & 0 & 0 & 2\cot(2r)
 \end{pmatrix}
\]
with respect to that decomposition. Therefore the tube of radius $0 < r <\pi/2$ around ${\mathbb C}P^k$ has four distinct constant principal curvatures $-\tan(r)$, $0$, $\cot(r)$ and $2\cot(2r)$ (unless $m=2$ in which case there are only two distinct constant principal curvatures $0$ and $2\cot(2r)$). The corresponding principal curvature spaces are $[AN]$, $T_{[z]}{\mathbb C}P^k \ominus [AN]$,  $\nu_{[z]}{\mathbb C}P^k \ominus [N]$ and ${\mathbb R}JN$ respectively, where we identify the subspaces obtained by parallel translation along $\gamma$ from $[z]$ to $\gamma(r)$.
Note that the parallel translate of $[AN]$ corresponds to ${\mathcal C} \ominus {\mathcal Q}$, the parallel translate of $[N]$ corresponds to ${\mathbb C}\nu M$, and the parallel translate of ${\mathbb R}JN$ corresponds to ${\mathcal F}$. Moreover, we have $A(T_{[z]}{\mathbb C}P^k \ominus [AN]) = \nu_{[z]}{\mathbb C}P^k \ominus [N]$.

Since $JN$ is a principal curvature vector we also conclude that every tube around ${\mathbb C}P^k$ is a Hopf hypersurface. We also see that all principal curvature spaces orthogonal to ${\mathbb R}JN$ are $J$-invariant. Thus, if $\phi$ denotes the structure tensor field on the tube which is induced by $J$, we get $S\phi = \phi S$. Since the complex structure on $Q^m$ is parallel we have 
\[
g(\nabla_X\xi,Y) + g(X,\nabla_Y\xi) = g((S\phi - \phi S)X,Y)
\]
for all $X,Y \in TM$. As $\xi$ is a Killing vector field if and only
if $\nabla\xi$ is a skew-symmetric tensor field, we see that the Reeb flow on $M$ is isometric if and only if $S\phi = \phi S$.

We summarize the previous discussion in the following proposition.

\begin{prop}
Let $M$ be the tube of radius $0 < r <\pi/2$ around the totally geodesic ${\mathbb C}P^k$ in $Q^{2k}$, $k \geq 2$. Then the following statements hold:
\begin{enumerate}
\item[1.] {$M$ is a Hopf hypersurface.}
\item[2.] {The tangent bundle $TM$ and the normal bundle $\nu M$ of $M$ consist of ${\mathfrak A}$-isotropic singular tangent vectors of $Q^{2k}$.}
\item[3.] {$M$ has four distinct constant principal curvatures. Their values and corresponding principal curvature spaces and multiplicities are given in the following table:
\medskip
\begin{center}
\begin{tabular}{|c|c|c|}
\hline
\mbox{principal curvature} & \mbox{eigenspace}  & \mbox{multiplicity}  \\
\hline
$2\cot(2r)$ & ${\mathcal F}$ & $1$ \\
$0$ & ${\mathcal C} \ominus {\mathcal Q}$ & $2$ \\
$-\tan(r)$ & $T{\mathbb C}P^k \ominus ({\mathcal C} \ominus {\mathcal Q})$ & $2k-2$\\
$\cot(r)$ & $\nu{\mathbb C}P^k \ominus {\mathbb C}\nu M$ & $2k-2$ \\
\hline
\end{tabular},
\end{center}
\medskip
}
\noindent The real structure $A$ determined by the ${\mathfrak A}$-isotropic unit normal vector at $[z]$ maps $T_{[z]}{\mathbb C}P^k \ominus ({\mathcal C}_{[z]} \ominus {\mathcal Q}_{[z]})$ onto $\nu_{[z]}{\mathbb C}P^k \ominus {\mathbb C}\nu_{[z]} M$, and vice versa.
\item[4.] {The shape operator $S$ of $M$ and the structure tensor field $\phi$ of $M$ commute with each other, that is, $S\phi = \phi S$.}
\item[5.] {The Reeb flow on $M$ is an isometric flow.}
\end{enumerate}
\end{prop}

\section {The Codazzi equation and some consequences}

\medskip
From the explicit expression of the Riemannian curvature tensor of the complex quadric $Q^m$ we can easily derive the Codazzi equation for a real hypersurface $M \subset Q^m$:
\begin{eqnarray*}
g((\nabla_XS)Y - (\nabla_YS)X,Z) & = & \eta(X)g(\phi Y,Z) - \eta(Y) g(\phi X,Z) - 2\eta(Z) g(\phi X,Y) \\
 & & + \rho(X)g(BY,Z) - \rho(Y)g(BX,Z)\\
& & - \eta(BX)g(BY,\phi Z) - \eta(BX)\rho(Y)\eta(Z)  \\
& & + \eta(BY)g(BX,\phi Z) + \eta(BY)\rho(X)\eta(Z)  .
\end{eqnarray*}
We now assume that $M$ is a Hopf hypersurface. Then we have 
\[
S\xi = \alpha \xi
\]
with the smooth function $\alpha = g(S\xi,\xi)$ on $M$.
Inserting $Z = \xi$ into the Codazzi equation leads to
\[
g((\nabla_XS)Y - (\nabla_YS)X,\xi) =  - 2 g(\phi X,Y)  + 2\rho(X)\eta(BY) - 2\rho(Y)\eta(BX).
\]
On the other hand, we have
\begin{eqnarray*}
 & & g((\nabla_XS)Y - (\nabla_YS)X,\xi) \\
& = & g((\nabla_XS)\xi,Y) - g((\nabla_YS)\xi,X) \\
& = & d\alpha(X)\eta(Y) - d\alpha(Y)\eta(X) + \alpha g((S\phi + \phi S)X,Y) - 2g(S \phi SX,Y).
\end{eqnarray*}
Comparing the previous two equations and putting $X = \xi$ yields
$$
d\alpha(Y)  =  d\alpha(\xi)\eta(Y)  - 2\delta\rho(Y).
$$
Reinserting this into the previous equation yields
\begin{eqnarray*}
 & & g((\nabla_XS)Y - (\nabla_YS)X,\xi) \\
& = &  2\delta\eta(X)\rho(Y) - 2\delta \rho(X)\eta(Y) + \alpha g((\phi S + S\phi)X,Y) - 2g(S \phi SX,Y) .
\end{eqnarray*}
Altogether this implies
\begin{eqnarray*}
0 & = & 2g(S \phi SX,Y) - \alpha g((\phi S + S\phi)X,Y) - 2 g(\phi X,Y) \\
& & + 2\delta \rho(X)\eta(Y) + 2\rho(X)\eta(BY) - 2\rho(Y)\eta(BX) - 2\delta\eta(X)\rho(Y) \\
& = & g(2S \phi S - \alpha(\phi S + S\phi) - 2\phi) X ,Y)\\
& & + 2\rho(X)\eta(BY+\delta Y) - 2\rho(Y)\eta(BX+\delta X) \\
& = & g(2S \phi S - \alpha(\phi S + S\phi) - 2\phi) X ,Y)\\
& & + 2\rho(X)g(Y,B\xi+\delta \xi) - 2g(X,B\xi+\delta \xi)\rho(Y).
\end{eqnarray*}
If $AN = N$ we have $\rho = 0$, otherwise we can use Lemma \ref{notprincipal} to calculate
\[
\rho(Y) = g(Y,AN) = g(Y,AJ\xi) = -g(Y,JA\xi) = -g(Y,JB\xi) = -g(Y,\phi B\xi).
\]
Thus we have proved

\begin{lm}\label{Cod1}
Let $M$ be a Hopf hypersurface in $Q^m$, $m \geq 3$. Then we have
\[
(2S \phi S - \alpha(\phi S + S\phi) - 2\phi) X = 
- 2\rho(X)(B\xi+\delta \xi) - 2g(X,B\xi+\delta \xi)\phi B\xi .
\]
\end{lm}

If $N$ is ${\mathfrak A}$-principal we can choose a real structure $A \in {\mathfrak A}$ such that $AN = N$. Then we have $\rho = 0$ and $\phi B\xi = -\phi \xi = 0$, and therefore
\[ 
2S \phi S - \alpha(\phi S + S\phi) = 0.
\]
If $N$ is not ${\mathfrak A}$-principal we can choose a real structure $A \in {\mathfrak A}$ as in Lemma \ref{notprincipal} and get
\begin{eqnarray*}
& & \rho(X)(B\xi+\delta \xi) + g(X,B\xi+\delta \xi)\phi B\xi \\
& = & -g(X,\phi(B\xi + \delta\xi))(B\xi+\delta \xi) + g(X,B\xi+\delta \xi)\phi (B\xi + \delta\xi)\\
& = & ||B\xi + \delta\xi||^2( g(X,U)\phi U -g(X,\phi U)U ) \\
& = & \sin^2(2t)( g(X,U)\phi U -g(X,\phi U)U ),
\end{eqnarray*}
which is equal to $0$ on ${\mathcal Q}$ and equal to $\sin^2(2t)\phi X$ on ${\mathcal C} \ominus {\mathcal Q}$. Altogether we have proved:

\begin{lm}\label{Cod2}
Let $M$ be a Hopf hypersurface in $Q^m$, $m \geq 3$. Then the tensor field
\[ 2S\phi S - \alpha (\phi S + S\phi) \]
leaves ${\mathcal Q}$ and ${\mathcal C} \ominus {\mathcal Q}$ invariant and we have
\[ 2S\phi S - \alpha (\phi S + S\phi) = 2\phi \ {\rm on}\  {\mathcal Q} \]
and
\[ 2S\phi S - \alpha (\phi S + S\phi) = 2\delta^2\phi \ {\rm on}\  {\mathcal C} \ominus {\mathcal Q}. \]
\end{lm}

We will now prove that the principal curvature $\alpha$ of a Hopf hypersurface is constant if the normal vectors are ${\mathfrak A}$-isotropic. Assume that $N$ is ${\mathfrak A}$-isotropic everywhere. Then we have $\delta = 0$ and we get
\[
Y\alpha  =  d\alpha(\xi)\eta(Y)
\]
for all $Y \in TM$.
Since ${\rm grad}^M \alpha =d\alpha(\xi)\xi$, we can compute the Hessian ${\rm hess}^M \alpha$ by
\[
({\rm hess}^M \alpha)(X,Y) = g(\nabla_X{\rm grad}^M \alpha , Y)
= d(d\alpha(\xi))(X)\eta(Y) +d\alpha(\xi) g(\phi SX,Y).
\]
As ${\rm hess}^M \alpha$ is a symmetric bilinear form, the previous equation implies
\[
d\alpha(\xi)g((S\phi + \phi S)X,Y) = 0
\]
for all vector fields $X,Y$ on $M$ which are tangential to ${\mathcal C}$.

Now let us assume that $S{\phi}+{\phi}S=0$. For every principal curvature vector $X \in {\mathcal C}$ such that $SX={\lambda}X$ this implies $S{\phi}X=-{\phi}SX=-{\lambda}{\phi}X$. We assume $||X|| = 1$ and put $Y = \phi X$. Using Lemma \ref{Cod1} we get
\[
1 \leq \lambda^2 + 1 = \rho(X)\eta(B\phi X) - \rho(\phi X)\eta(BX) = g(X,AN)^2 + g(X,A\xi)^2 =  ||X_{{\mathcal C} \ominus {\mathcal Q}}||^2 \leq 1,
\]
where $X_{{\mathcal C} \ominus {\mathcal Q}}$ denotes the orthogonal projection of $X$ onto ${\mathcal C} \ominus {\mathcal Q}$. This
implies $||X_{{\mathcal C} \ominus {\mathcal Q}}||^2 = 1$
for all principal curvature vectors $X \in {\mathcal C}$ with $||X|| = 1$. This is only possible if ${\mathcal C} = {\mathcal C} \ominus {\mathcal Q}$, or equivalently, if ${\mathcal Q} = 0$. Since $m \geq 3$ this is not possible. Hence we must have $S{\phi}+{\phi}S \neq 0$ everywhere and therefore $d\alpha(\xi) = 0$, which implies ${\rm grad}^M \alpha = 0$. Since $M$ is connected this implies that $\alpha$ is constant. Thus we have proved:

\medskip

\begin{lm} \label{alphaconstant}
Let $M$ be a Hopf hypersurface in $Q^m$, $m \geq 3$, such that the normal vector
field $N$ is ${\mathfrak A}$-isotropic everywhere. Then $\alpha$ is
constant. 
\end{lm}

\section{Proof of Theorem \ref{mainresult}}

Let $M$ be a real hypersurface in $Q^m$, $m \geq 3$, with isometric Reeb flow. As we have seen above, this geometric condition is equivalent to the algebraic condition  $S\phi = \phi S$. Applying this equation to $\xi$ gives $0 = S\phi \xi = \phi S\xi$, which implies that $S\xi = \alpha \xi$ with $\alpha = g(S\xi,\xi)$. Therefore any real hypersurface in $Q^m$ with isometric Reeb flow is a Hopf hypersurface.

Differentiating the equation $S\phi - \phi S = 0$ gives
\begin{eqnarray*}
0 & = & (\nabla_XS)\phi Y + S(\nabla_X\phi)Y - (\nabla_X\phi)SY - \phi(\nabla_XS)Y \\
& = & (\nabla_XS)\phi Y + S(\eta(Y)SX-g(SX,Y)\xi) - (\eta(SY)SX-g(SX,SY)\xi) - \phi(\nabla_XS)Y \\
& = & (\nabla_XS)\phi Y + \eta(Y)S^2X-\alpha g(SX,Y)\xi - \eta(SY)SX + g(SX,SY)\xi - \phi(\nabla_XS)Y. 
\end{eqnarray*}
If we define
\[
\Psi(X,Y,Z) = g((\nabla_XS)Y,\phi Z) + g((\nabla_XS)Z,\phi Y),
\]
the previous equation implies 
\[
\Psi(X,Y,Z) = \alpha\eta(Z)g(SX,Y) - \eta(Z)g(SX,SY) + \eta(SY)g(SX,Z) - \eta(Y)g(S^2X,Z).
\]
Evaluating $\Psi(X,Y,Z) + \Psi(Y,Z,X) - \Psi(Z,X,Y)$  leads to
\begin{eqnarray*}
2g((\nabla_XS)Y,\phi Z) & = & \Omega(X,Y,Z) - \Omega(Y,Z,X) + \Omega(Z,X,Y) \\
& & + 2\alpha\eta(Z)g(SX,Y) - 2\eta(Z)g(S^2X,Y),
\end{eqnarray*}
where
\[
\Omega(X,Y,Z) = g((\nabla_XS)Y - (\nabla_YS)X,\phi Z).
\]
The three $\Omega$-terms can be evaluated using the Codazzi equation, which leads to
\begin{eqnarray*}
& & 2g((\nabla_XS)Y,\phi Z) \\ 
& = &  \rho(X)(g(AY,\phi Z) - g(AZ,\phi Y)) + \eta(BX)(g(JAY,\phi Z) - g(JAZ,\phi Y))\\
& & - \rho(Y)(g(AX,\phi Z) + g(AZ,\phi X)) -  \eta(BY)(g(JAX,\phi Z) + g(JAZ,\phi X))\\
& & + \rho(Z)(g(AX,\phi Y) + g(AY,\phi X)) + \eta(BZ)(g(JAX,\phi Y) + g(JAY,\phi X))\\
& & + 2\eta(Z)g(\phi X,\phi Y) - 2\eta(Y)g(\phi X,\phi Z) + 2\alpha\eta(Z)g(SX,Y) - 2\eta(Z)g(S^2X,Y).
\end{eqnarray*}
Replacing $\phi Z$ by $JZ - \eta(Z)N$, and similarly for $X$ and $Y$, one can easily calculate that
\begin{eqnarray*}
g(JAY,\phi Z) - g(JAZ,\phi Y) & = & \eta(Y)\eta(BZ) - \eta(Z)\eta(BY), \\
g(JAX,\phi Z) + g(JAZ,\phi X) & = & 2g(BX,Z) - \eta(X)\eta(BZ) - \eta(Z)\eta(BX),\\
g(JAX,\phi Y) + g(JAY,\phi X) & = & 2g(BX,Y) - \eta(X)\eta(BY) - \eta(Y)\eta(BX).
\end{eqnarray*}
Inserting this into the previous equation gives
\begin{eqnarray*}
& & 2g((\nabla_XS)Y,\phi Z) \\ 
& = &  \rho(X)(g(BY,\phi Z) - g(BZ,\phi Y)) - \rho(Y)(g(BX,\phi Z) + g(BZ,\phi X)) \\
& & + \rho(Z)(g(BX,\phi Y) + g(BY,\phi X)) - 2\eta(BY)g(BX,Z) + 2\eta(BZ)g(BX,Y)  \\
& & + 2\eta(Z)g(\phi X,\phi Y) - 2\eta(Y)g(\phi X,\phi Z) + 2\alpha\eta(Z)g(SX,Y) - 2\eta(Z)g(S^2X,Y).
\end{eqnarray*}
Since
\begin{eqnarray*}
g(BY,\phi Z) - g(BZ,\phi Y)  & =  & \eta(Y)\rho(Z) - \eta(Z)\rho(Y), \\
g(BX,\phi Y) + g(BY,\phi X)  & =  & 2g(BX,\phi Y) + \eta(Y)\rho(X) - \eta(X)\rho(Y), \\
g(BX,\phi Z) + g(BZ,\phi X)  & =  & 2g(BX,\phi Z) + \eta(Z)\rho(X) - \eta(X)\rho(Z),
\end{eqnarray*}
we get
\begin{eqnarray*}
& & g((\nabla_XS)Y,\phi Z) \\ 
& = &  \rho(X)\eta(Y)\rho(Z) - \rho(X)\rho(Y)\eta(Z) - \rho(Y)g(BX,\phi Z)  + \rho(Z)g(BX,\phi Y) \\
& &  - \eta(BY)g(BX,Z) + \eta(BZ)g(BX,Y)  \\
& & + \eta(Z)g(\phi X,\phi Y) - \eta(Y)g(\phi X,\phi Z) + \alpha\eta(Z)g(SX,Y) - \eta(Z)g(S^2X,Y).
\end{eqnarray*}
Replacing $Z$ by $\phi Z$ and using $\phi^2 Z = - Z + \eta(Z)\xi$ gives
\begin{eqnarray*}
& & -g((\nabla_XS)Y,Z) + \eta(Z)g((\nabla_XS)Y,\xi) \\ 
& = &  \rho(X)\eta(Y)\rho(\phi Z)  + \rho(Y)g(BX,Z) - \rho(Y)\eta(Z)g(BX,\xi) + \rho(\phi Z)g(BX,\phi Y) \\
& &  - \eta(BY)g(BX,\phi Z) + \eta(B\phi Z)g(BX,Y) + \eta(Y)g(\phi X,Z).
\end{eqnarray*}
Since
\[
g((\nabla_XS)Y,\xi) = d\alpha(X)\eta(Y) + \alpha g(S \phi X,Y) - g(S^2 \phi X,Y),
\]
this implies
\begin{eqnarray*}
& & g((\nabla_XS)Y,Z)  \\ 
& = &  d\alpha(X)\eta(Y)\eta(Z) + \alpha \eta(Z) g(S \phi X,Y) - \eta(Z) g(S^2 \phi X,Y) \\
& &  - \rho(X)\eta(Y)\rho(\phi Z)  - \rho(Y)g(BX,Z) + \rho(Y)\eta(Z)g(BX,\xi) - \rho(\phi Z)g(BX,\phi Y) \\
& &  + \eta(BY)g(BX,\phi Z) - \eta(B\phi Z)g(BX,Y) - \eta(Y)g(\phi X,Z).
\end{eqnarray*}
From this we get an explicit expression for the covariant derivaive of the shape operator,
\begin{eqnarray*}
(\nabla_XS)Y  & = &
\{d\alpha(X)\eta(Y) + \alpha g(S \phi X,Y) -  g(S^2 \phi X,Y)  + \delta\eta(Y)\rho(X)  \\
& & \quad + \delta g(BX,\phi Y) + \eta(BX)\rho(Y)\} \xi\\
& &  + \{\eta(Y)\rho(X) +g(BX,\phi Y)\} B\xi + g(BX,Y)\phi B\xi  \\
& & - \rho(Y)BX - \eta(Y)\phi X - \eta(BY)\phi BX  .
\end{eqnarray*}
Putting $Y = \xi$ and $X \in {\mathcal C}$ then leads to
\begin{eqnarray*}
\alpha S \phi X - S^2 \phi  X & = &  \delta \rho(X) \xi + \rho(X)  B\xi + \eta(BX)\phi B\xi - \phi X + \delta\phi BX  .
\end{eqnarray*}
On the other hand, from Lemma \ref{Cod1} we get
\begin{eqnarray*}
\alpha S \phi X - S^2 \phi  X & = &  \delta \rho(X) \xi + \rho(X)  B\xi + \eta(BX)\phi B\xi - \phi X   .
\end{eqnarray*}
Comparing the previous two equations leads to
\[
\delta\phi BX = 0
\]
for all $X \in {\mathcal C}$. Let us first assume that $\delta \neq 0$. Then we have $\phi BX = 0$ for all $X \in {\mathcal C}$, which implies $BX = \eta(BX)\xi$ for all $X \in {\mathcal C}$, and therefore
\[
AX = BX + \rho(X)N = \eta(BX)\xi + \rho(X)N
\]
for all $X \in {\mathcal C}$. This implies $A({\mathcal C}) \subset [N]$, which gives a contradiction since $A$ is an isomorphism everywhere and the rank of ${\mathcal C}$ is equal to $2(m-1)$ and $m \geq 3$. Therefore we must have $\delta = 0$, which means that $N$ is ${\mathfrak A}$-isotropic. We thus have proved the following proposition.

\begin{prop} \label{Aisotropic}
Let $M$ be a real hypersurface in $Q^m$, $m \geq 3$, with isometric Reeb flow. Then the normal vector
field $N$ is ${\mathfrak A}$-isotropic everywhere. 
\end{prop}

From Proposition \ref{Aisotropic} and Lemma \ref{alphaconstant} we conclude that the principal curvature function $\alpha$ is constant.
Since the normal vector field $N$ is ${\mathfrak A}$-isotropic, the vector fields $N$, $\xi$, $B\xi$ and $\phi B\xi$ are pairwise orthonormal. This implies that
\[
{\mathcal C} \ominus {\mathcal Q} = [B\xi].
\]
From Lemma \ref{Cod2} we know that the tensor field $2S\phi S - \alpha (\phi S + S\phi) = 2\phi (S^2 - \alpha S)$
leaves ${\mathcal Q}$ and ${\mathcal C} \ominus {\mathcal Q}$ invariant and 
\[ \phi (S^2 - \alpha S) = \phi \   {\rm on}\  {\mathcal Q} \ \ {\rm and} \ \ 
\phi (S^2 - \alpha S) = 0 \ {\rm on}\  {\mathcal C} \ominus {\mathcal Q}. \]
Since $\phi$ is an isomorphism of ${\mathcal Q}$ and of ${\mathcal C} \ominus {\mathcal Q}$ this implies 
\[ S^2 - \alpha S = I \ {\rm on}\  {\mathcal Q} \ \ {\rm and}\ \ 
S^2 - \alpha S = 0 \ {\rm on}\  {\mathcal C} \ominus {\mathcal Q}. \]
As $M$ is a Hopf hypersurface we have $S({\mathcal C}) \subset {\mathcal C}$. Let $X \in {\mathcal C}$ be a principal curvature vector of $M$ with corresponding principal curvature $\lambda$, that is, $SX = \lambda X$. We decompose $X$ into $X = Y + Z$ with $Y \in {\mathcal Q}$ and $Z \in {\mathcal C} \ominus {\mathcal Q}$. Then we get
\[
(\lambda^2 -\alpha\lambda)Y + (\lambda^2 -\alpha\lambda)Z = (\lambda^2 -\alpha\lambda)X = (S^2 - \alpha S)X = Y.
\]
If $\lambda^2 - \alpha\lambda = 0$ we must have $Y = 0$ and therefore $X \in {\mathcal C} \ominus {\mathcal Q}$. If $\lambda^2 - \alpha\lambda \neq 0$ we must have $Z = 0$ and therefore $X \in {\mathcal Q}$. Altogether this implies

\begin{prop} \label{Qinvariance}
Let $M$ be a real hypersurface in $Q^m$, $m \geq 3$, with isometric Reeb flow. Then the distributions ${\mathcal Q}$ and ${\mathcal C} \ominus {\mathcal Q} = [B\xi]$ are invariant under the shape operator $S$ of $M$.
\end{prop}

Assume that $SX = \lambda X$ with $X \in {\mathcal C}$. From Lemma \ref{Cod2} we get $\lambda^2 - \alpha\lambda - 1 = 0$ if $X \in {\mathcal Q}$ and $\lambda^2 - \alpha\lambda = 0$ if $X \in {\mathcal C} \ominus {\mathcal Q}$. Recall that $\alpha$ is constant. We put $\alpha = 2\cot(2r)$ with $0 < r < \frac{\pi}{2}$ and define $T_\alpha = {\mathbb R}\xi = {\mathcal F}$. Then the solutions of $x^2 -\alpha x - 1 = 0$ are $\lambda = \cot(r)$ and $\mu = -\tan(r)$. We denote by $T_\lambda$ and $T_\mu$ the subbundles of ${\mathcal Q}$ consisting of the corresponding principal curvature vectors. The rank of ${\mathcal C} \ominus {\mathcal Q}$ is equal to $2$ and ${\mathcal C} \ominus {\mathcal Q}$ is both $S$- and $\phi$-invariant. Therefore, since $S\phi = \phi S$, there is exactly one principal curvature $\nu$ on ${\mathcal C} \ominus {\mathcal Q}$ which is equal to either $\nu = 0$ or $\nu = \alpha$. We define $T_\nu = {\mathcal C} \ominus {\mathcal Q}$. Note that, since $S\phi = S\phi$, we have $JT_\rho = \phi T_\rho = T_\rho$ for $\rho \in \{\lambda,\mu,\nu\}$. 

\medskip
According to Lemma \ref{Cod1} we have
\[
(\alpha S\phi - S^2 \phi)X =  
-\phi X  + \rho(X)B\xi + g(X,B\xi)\phi B\xi = 
- \phi X  + \eta(B\phi X)B\xi + \eta(BX)\phi B\xi.
\]
Inserting this into the above expression for the covariant derivative of $S$ and replacing $\rho(X)$ and $\rho(Y)$ by $\eta(B\phi X)$ and $\eta(B \phi Y)$ respectively leads to
\begin{eqnarray*}
(\nabla_XS)Y  & = &
\{g((\alpha S \phi - S^2 \phi) X,Y)  + \eta(BX)\eta(B\phi Y)\} \xi\\
& &  + \{\eta(Y)\eta(B \phi X) +g(BX,\phi Y)\} B\xi + g(BX,Y)\phi B\xi  \\
& & - \eta(B \phi Y)BX - \eta(Y)\phi X - \eta(BY)\phi BX  \\
& = & \{\eta(B\phi X)\eta(BY) - g(\phi X,Y)\} \xi\\
& &  + \{\eta(B\phi X)\eta(Y)-g(\phi BX,Y)\} B\xi + g(BX,Y)\phi B\xi  \\
& & - \eta(B\phi Y)BX - \eta(Y)\phi X - \eta(BY)\phi BX  .
\end{eqnarray*}
This implies
\[ g(BX,Y) = g((\nabla_XS)Y,\phi B\xi).\]
Let $X$ and $Y$ be sections in $T_\lambda$. Then the previous equation implies
\[
g(BX,Y) = (\lambda - \nu)g(\nabla_XY,\phi B\xi) = (\nu - \lambda)g(\nabla_X(\phi B\xi),Y) = (\nu - \lambda)g(\phi\nabla_X(B\xi),Y).
\]
Since $B\xi = A\xi$, the Gauss formula for $M$ in $Q^m$ gives
\[
\nabla_X(B\xi) = \nabla_X(A\xi) = \bar\nabla_X(A\xi),
\]
where $\bar\nabla$ is the Riemannian connection of $Q^m$. Since ${\mathfrak A}$ is a parallel subbundle of ${\rm End}(TQ^m)$, there exists $A^\prime \in {\mathfrak A}$ such that $\bar\nabla_XA = A^\prime$, and we get
\[
\bar\nabla_X(A\xi) = A\bar\nabla_X\xi + A^\prime\xi = A\nabla_X\xi + A^\prime\xi = A\phi SX + A^\prime \xi = \lambda A\phi X + A^\prime\xi.
\]
Since $A^\prime$ differs from $A$ by a complex scalar we have $A^\prime \xi \in [A\xi] = [B\xi]$ and thus $g(\phi A^\prime \xi, Y) = 0$.
Altogether we therefore get
\[
g(BX,Y) = (\nu - \lambda)\lambda g(\phi A\phi X,Y) = (\nu - \lambda)\lambda g(AX,Y) = (\nu - \lambda)\lambda g(BX,Y).
\]
Assume that $g(BX,Y) \neq 0$. Then we have $(\nu - \lambda)\lambda = 1$. Recall that $\nu \in \{0,\alpha\}$. If $\nu = 0$ we get $-\lambda^2 = 1$, and if $\nu = \alpha$ we get $-1 = (\alpha - \lambda)\lambda = 1$, both of which gives a contradiction. Therefore we must have $g(BX,Y) = 0$ for all sections $X$ and $Y$ in $T_\lambda$. The same argument can be repeated for $T_\mu$. Since ${\mathcal Q} = T_\lambda \oplus T_\mu$ and $B({\mathcal Q}) = {\mathcal Q}$ we conclude
\[
B(T_\lambda) = T_\mu\ \ {\rm and}\ \ B(T_\mu) = T_\lambda.
\]
Since $B = A$ on $T_\lambda$ and $T_\mu$ we can replace $B$ by $A$ here. As both $T_\lambda$ and $T_\mu$ are complex we see that ${\mathcal Q}$ and hence $Q^m$ must have even complex dimension. We summarize this in

\begin{prop} \label{Qinvariance}
Let $M$ be a real hypersurface in $Q^m$, $m \geq 3$, with isometric Reeb flow. Then $m$ is even, say $m = 2k$, and the real structure $A$ maps $T_\lambda$ onto $T_\mu$, and vice versa.
\end{prop}

For each point $[z] \in M$ we denote by $\gamma_{[z]}$ the geodesic in $Q^{2k}$ with $\gamma_{[z]}(0)=[z]$ and $\dot{\gamma}_{[z]}(0)=N_{[z]}$ and by $F$ the smooth map
\[
F: M \longrightarrow Q^m, [z] \longrightarrow \gamma_{[z]}(r).
\]
Geometrically, $F$ is the displacement of $M$ at distance $r$ in the direction of the normal vector field $N$. For each point $[z] \in M$ the differential $d_{[z]}F$ of $F$ at $[z]$ can be computed by using Jacobi vector fields as
\[
d_{[z]}F(X)=Z_X(r),
\]
where $Z_X$ is the Jacobi vector field along $\gamma_{[z]}$ with initial values $Z_X(0)=X$ and $Z_X'(0)=-SX$. Using the fact that $N$ is ${\mathfrak A}$-isotropic, we can calculate the normal Jacobi operator $R_N$ from the explicit expression of the curvature tensor of $Q^{2k}$:
\[
R_N Z = R(Z,N)N  =   Z - g(Z,N)N + 3g(Z,\xi)\xi   - g(Z,AN)AN - g(Z,A\xi)A\xi.
\]
It follows that $R_N$ has the three constant eigenvalues $0,1,4$ with corresponding eigenbundles $\nu M \oplus ({\mathcal C} \ominus {\mathcal Q}) = \nu M \oplus T_\nu$, ${\mathcal Q} = T_\lambda \oplus T_\mu$ and ${\mathcal F} = T_\alpha$.
This leads to the following expressions for the Jacobi vector fields along ${\gamma}_{[z]}$:
\begin{equation*}
Z_X(r)=\begin{cases}
(\cos(2r) - \frac{\alpha}{2}\sin(2r))E_X(r) & \text{if $X \in T_\alpha$,}\\
(\cos(r) - \rho\sin(r))E_X(r) & \text{if $X\in T_\rho$, $\rho\in\{\lambda,\mu\}$,}\\
(1 - \nu r)E_X(r) & \text{if $X\in T_\nu$},
\end{cases}
\end{equation*}
where $E_X$ is the parallel vector field along ${\gamma}_{[z]}$ with $E_X(0)=X$. 
This shows that ${\rm Ker}(dF) = T_\alpha \oplus T_\lambda$ and thus $F$ has constant rank $2k$,
where we take into account that in the case $\nu = \alpha$ the function $1 - 2r\cot(2r)$ has positive values for $0 < r < \frac{\pi}{2}$.
So, locally, $F$ is a submersion onto a submanifold $P$ of $Q^{2k}$ of real dimension $2k$. Moreover, the tangent space $T_{F([z])}P$ of $P$ at $F([z])$ is obtained by parallel translation of $(T_\mu \oplus T_\nu)([z])$ along $\gamma_{[z]}$. Since $T_\mu$ and $T_\nu$ are both $J$-invariant and $J$ is invariant under parallel translation along geodesics, the submanifold $P$ is a complex submanifold of $Q^{2k}$ of complex dimension $k$.

The vector $\eta_{[z]} = \dot{\gamma}_{[z]}(r)$ is a unit normal vector of $P$
at $F([z])$ and the shape operator $S_{\eta_{[z]}}$ of $P$ with respect to
$\eta_{[z]}$ can be calculated from the equation
\[
S_{\eta_{[z]}}Z_X(r) = - Z_X^\prime(r),
\]
where $X \in (T_\mu \oplus T_\nu)([z])$. The above expression for the Jacobi vector fields $Z_X$ implies $Z_X^\prime(r) = 0$ for $X \in T_\mu([z])$ and $X \in T_{\nu=0}([z])$, and therefore $S_{\eta_{[z]}} = 0$ if $\nu = 0$. If $\nu = \alpha$ we have $S_{\eta_{[z]}}E_X(r) = \frac{2\cot(2r)}{2r\cot(2r)-1} E_X(r)$. Since every complex submanifold of a K\"{a}hler manifold is minimal we must have $r = \frac{\pi}{4}$ and therefore $\nu = \alpha = 2\cot(2r) = 0$, which takes us back to the first case $\nu = 0$.
The vectors of the form $\eta_{[q]}$, $[q] \in F^{-1}(\{[z]\})$,
form an open subset of the unit
sphere in the normal space of $P$ at $F([z])$. Since $B_{\eta_{[q]}}$ vanishes
for all $\eta_{[q]}$ it follows that $P$ is a $k$-dimensional totally geodesic complex submanifold of $Q^{2k}$.
Rigidity of totally geodesic submanifolds now
implies that the entire submanifold
$M$ is an open part of a tube of radius $r$
around a $k$-dimensional connected, complete,
totally geodesic complex submanifold $P$ of $Q^{2k}$.

Klein classified in \cite{K} the totally geodesic submanifolds in complex quadrics. According to his classification the submanifold $P$ is either a totally geodesic $Q^k \subset Q^{2k}$ or a totally geodesic ${\mathbb C}P^k \subset Q^{2k}$. The normal spaces of $Q^k$ are Lie triple systems and the corresponding totally geodesic submanifolds of $Q^{2k}$ are again $k$-dimensional quadrics. Since $k \geq 2$ it follows that the normal spaces of $Q^k$ contain all types of tangent vectors of $Q^{2k}$. This implies that the normal bundle of the tubes around $Q^k$ contains regular and singular tangent vectors of $Q^{2k}$. Since the normal bundle of $M$ consists of ${\mathfrak A}$-isotropic tangent vectors only we conclude that $P$ is congruent to ${\mathbb C}P^k$. It follows that $M$ is congruent to an open part of a tube around ${\mathbb C}P^k$. This concludes the proof of Theorem \ref{mainresult}.

\end{document}